\documentclass[12pt, a4paper]{amsart}
\usepackage[utf8]{inputenc}
\usepackage{amsthm}
\usepackage{mathrsfs}  
\usepackage[none]{hyphenat}
\usepackage[english]{babel}
\usepackage{setspace}
\usepackage{pinlabel}

\usepackage[a4paper, top=3cm, bottom=4cm,left=3cm,right=3cm]{geometry}

\usepackage{comment}

\usepackage{pgf,tikz}
\usepackage{microtype}

\usepackage{tikz-cd}
\usepackage{tikz,pgfplots}
\usetikzlibrary{positioning, angles}
\usetikzlibrary{arrows}
\usetikzlibrary{decorations.pathmorphing, decorations.pathreplacing}
\usetikzlibrary{decorations.markings}
\usetikzlibrary{calc,shapes}
\usepackage{mathtools}
\usepackage{lipsum}
\usepackage{comment}
\usepackage{float}    
\usetikzlibrary{shapes.geometric}
\usepackage{tikz}
\usepackage{graphicx}
\let\phi\varphi

\usepackage{amsmath, amsthm, amssymb, amsfonts, xfrac}
\newtheorem{prop}{\textsc{Proposition}}[section]

\newtheorem{thm}[prop]{\textsc{Theorem}}

\theoremstyle{remark}

\newcommand{\Z}{\mathbb{Z}}

\newcommand{\C}{\widetilde{\mathcal{C}}}

\DeclareMathOperator{\Fix}{Fix}

\makeatletter
\def\blfootnote{\xdef\@thefnmark{}\@footnotetext}
\makeatother

\usepackage{url}
\usepackage{hyperref}
\title{The equivariant concordance group is not abelian}
\author{Alessio Di Prisa}
\address{Dipartimento di Matematica, Pisa, Italy \vskip.05in}
\email{a.diprisa@studenti.unipi.it}

\begin{document}

\maketitle
\begin{abstract}
We prove that the equivariant concordance group $\C$ is not abelian by exhibiting an infinite family of nontrivial commutators.
\end{abstract}
\section{Introduction}\label{intro}
A knot $K\subset S^3$ is said to be \emph{strongly invertible} if there is an orientation preserving smooth involution $\rho$ of $S^3$ such that $\rho(K)=K$ and $\rho$ reverses the orientation on $K$.
By introducing the notion of a \emph{direction} on a strongly invertible knot, Sakuma \cite{sakuma} was able to define unambiguously an operation of \emph{equivariant connected sum} between \emph{directed strongly invertible knots}. Moreover, he proved that this operation is not commutative, which is in stark contrast with the usual connected sum of oriented knots.
In the same paper, the author defined the \emph{equivariant concordance group} $\C$ as the quotient of the set of directed strongly invertible knots under the equivalence relation of \emph{equivariant concordance}, with the group operation naturally induced by the equivariant connected sum.

Since the equivariant connected sum is not commutative, it is a natural question whether the equivariant concordance group $\C$ is abelian.
Alfieri and Boyle \cite{alfieri_boyle} speculate that $\C$ contains a copy of $\Z*\Z$, and hence is nonabelian.
Recently Dai, Mallick and Stoffregen \cite{dai_mallick_stoffregen} define a homomorphism
$$h_{\tau,\iota}:\C\longrightarrow\mathfrak{K}_{\tau,\iota},$$
where $\mathfrak{K}_{\tau,\iota}$ is a (potentially) nonabelian group defined using the action induced by the strong involution of a knot $K$ on the knot Floer complex $\mathcal{CFK}(K)$.
Since the group $\mathfrak{K}_{\tau,\iota}$ is not (a priori) commutative, the authors observe that $h_{\tau,\iota}$ could in principle lead to a negative answer to the open question of whether $\C$ is abelian.

In this paper we present a family of examples which answer the question in the negative, showing that the equivariant concordance group $\C$ is indeed not abelian.
First of all, we recall some definitions and facts (see \cite{boyle_issa} for the details).
Given a knot $K$ with a strong inversion $\rho$ we say that $\Fix(\rho)$ is the \emph{axis} of the inversion. The axis is an unknotted $S^1$ which meets $K$ in two points.
A \emph{directed strongly invertible knot} $K$ is a knot with a strong inversion $\rho$ and the additional data of the choice of one of the components of $\Fix(\rho)\setminus K$ (a \emph{half-axis}) and an orientation on $\Fix(\rho)$.
Note that the points in $K\cap\Fix(\rho)$ have a natural order: the initial point of the chosen half-axis is the first fixed point while the other end is the second fixed point.
Two directed strongly inverible knots $K$ and $J$, with strong inversions $\rho_K$ and $\rho_J$, are \emph{equivariantly isotopic} if there exists an orientation preserving diffeomorphism $\phi$ of $S^3$ such that $\phi(K)=J$, $\phi\circ\rho_K=\rho_J\circ\phi$ and $\phi$ preserves the chosen oriented half-axis.
Given two directed strongly invertible knots $K$ and $J$, we denote their equivariant connected sum by $K\widetilde{\#}J$. The directed strongly invertible knot $K\widetilde{\#}J$ is obtained by cutting $K$ at its second fixed point and $J$ at its first fixed point, gluing the two knots and axes equivariantly, in a way that is compatible with the orientations on the axes, and choosing the half-axis of $K\widetilde{\#}J$ to be the union of the half-axes of $K$ and $J$.
The inverse $K^{-1}$ of a directed strongly invertible knot $K$ in $\C$ is represented by the mirror of $K$ with the same strong inversion and chosen half-axis, but with the opposite orientation on the axis of the strong inversion.

\begin{figure}[ht]
\centering
\begin{tikzpicture}

\node[anchor=south west,inner sep=0] at (0,0){\includegraphics[scale=0.3]{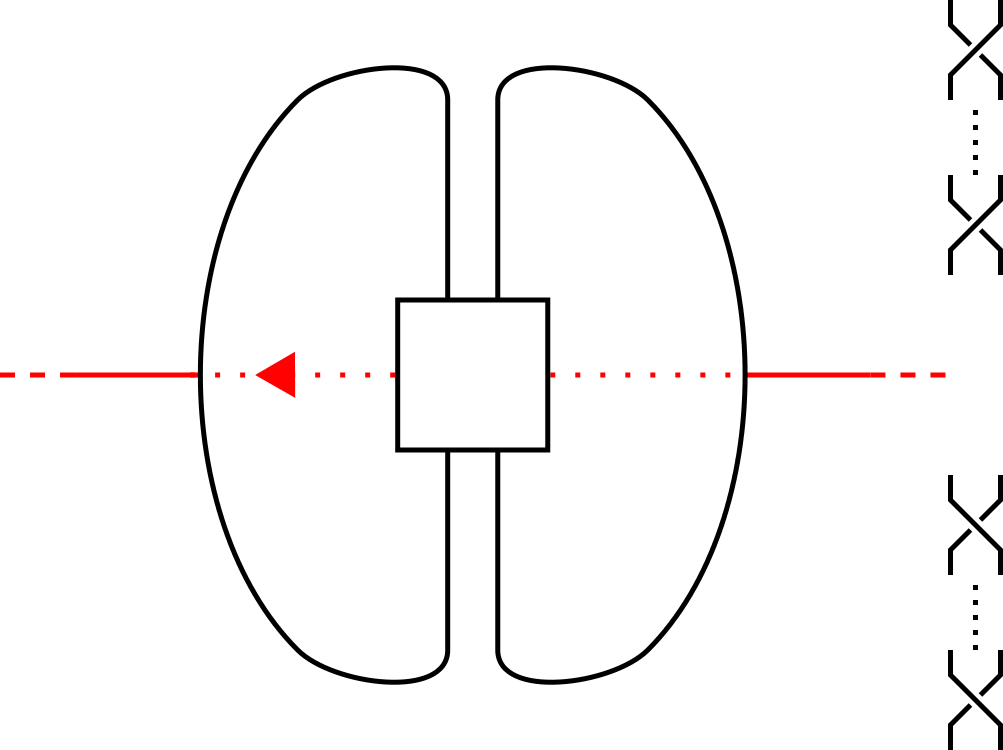}};
\node[label={$p$}] at (3.7,2.5){};
\node[label={$p>0$}] at (8.6,4.4){};
\node[label={$p<0$}] at (8.6,0.6){};

\end{tikzpicture}
    \caption{The directed strongly invertible knot $K_p$. The strong inversion is the $\pi$-rotation around the horizontal axis (colored red).}
    \label{torusknot}
\end{figure}
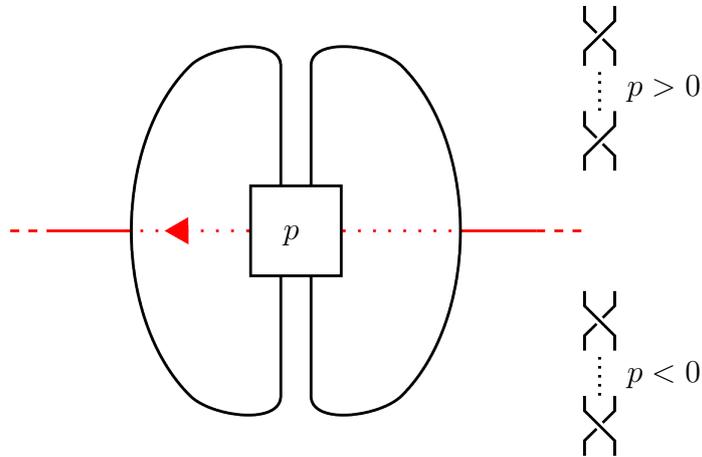

Let now $p$ be an odd integer and define $K_p$ to be the directed strongly invertible knot given by the torus knot $T_{2,p}$ with the orientation on the axis of the strong inversion described in Figure \ref{torusknot} and chosen half-axis given by the solid one in the figure. Hence first fixed point of $K_p$ is the one on the left in Figure \ref{torusknot}, while the second one is on the right.

\begin{thm}\label{teorema}
Let $p$ and $q$ be odd integers such that $1<p<q$. Then the commutator $[K_p,K_q]$ is not equivariantly slice. In particular the equivariant concordance group $\C$ is not abelian.
\end{thm}

\section{Proof of Theorem \ref{teorema}}\label{proof}
\begin{figure}[ht]
\centering
\begin{tikzpicture}

\node[anchor=south west,inner sep=0] at (0,0){\includegraphics[scale=0.3]{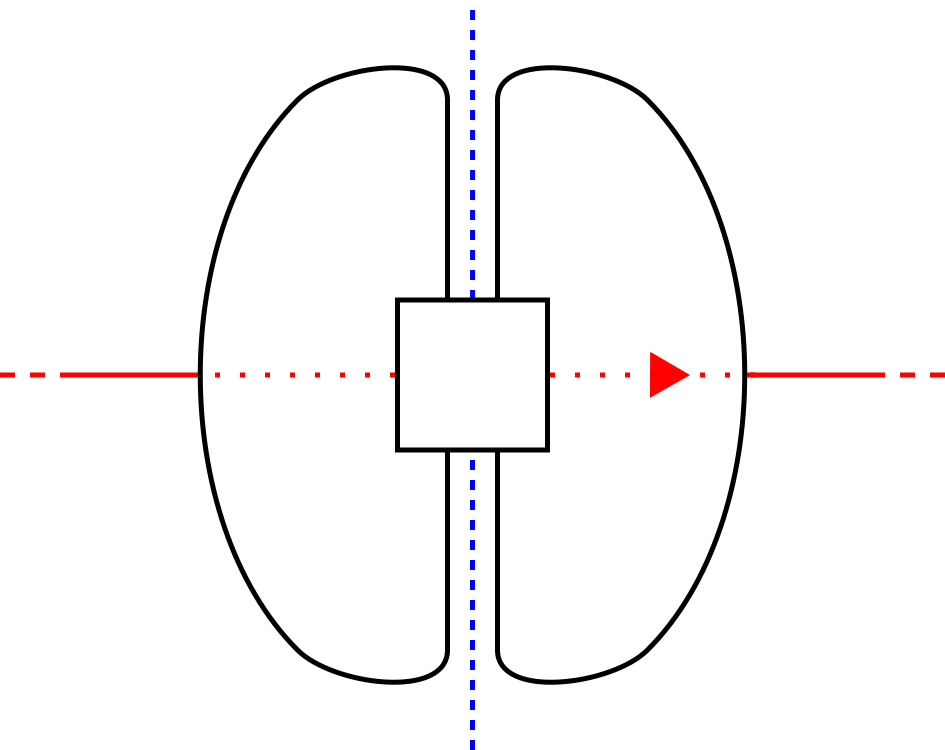}};
\node[label={$-p$}] at (3.7,2.5){};

\end{tikzpicture}
    \caption{The directed strongly invertible knot $K_p^{-1}$. The $\pi$-rotation around the vertical axis (colored blue) shows that $K_p^{-1}=K_{-p}$.}
    \label{inverse}
\end{figure}

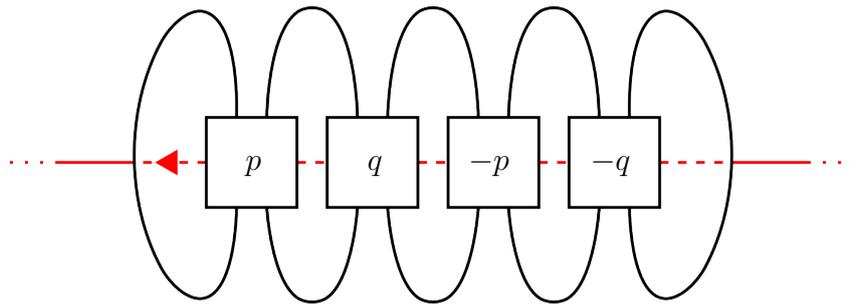
\begin{figure}[ht]
\centering
\begin{tikzpicture}

\node[anchor=south west,inner sep=0] at (0,0){\includegraphics[scale=0.3]{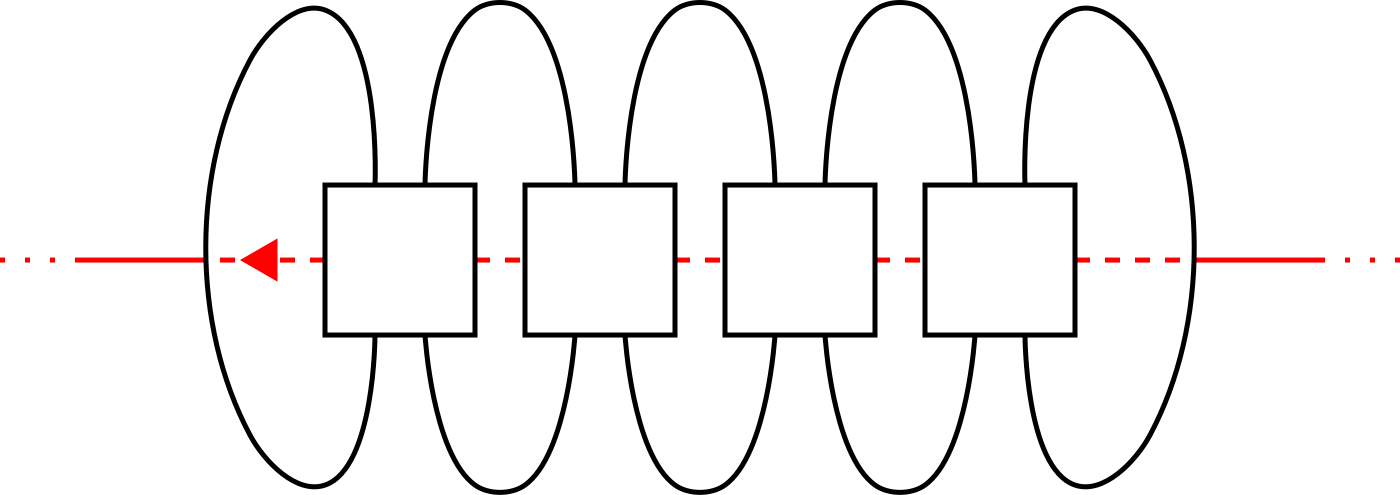}};
\node[label={$p$}] at (3.2,1.4){};
\node[label={$q$}] at (4.8,1.4){};
\node[label={$-p$}] at (6.3,1.4){};
\node[label={$-q$}] at (7.9,1.4){};

\end{tikzpicture}
    \caption{The commutator $[K_p,K_q]$.}
    \label{commutator}
\end{figure}

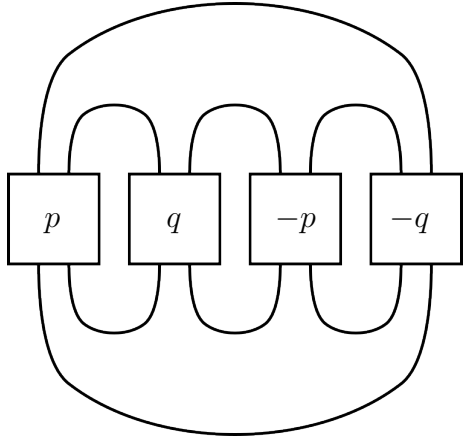
\begin{figure}[ht]
\centering
\begin{tikzpicture}

\node[anchor=south west,inner sep=0] at (0,0){\includegraphics[scale=0.3]{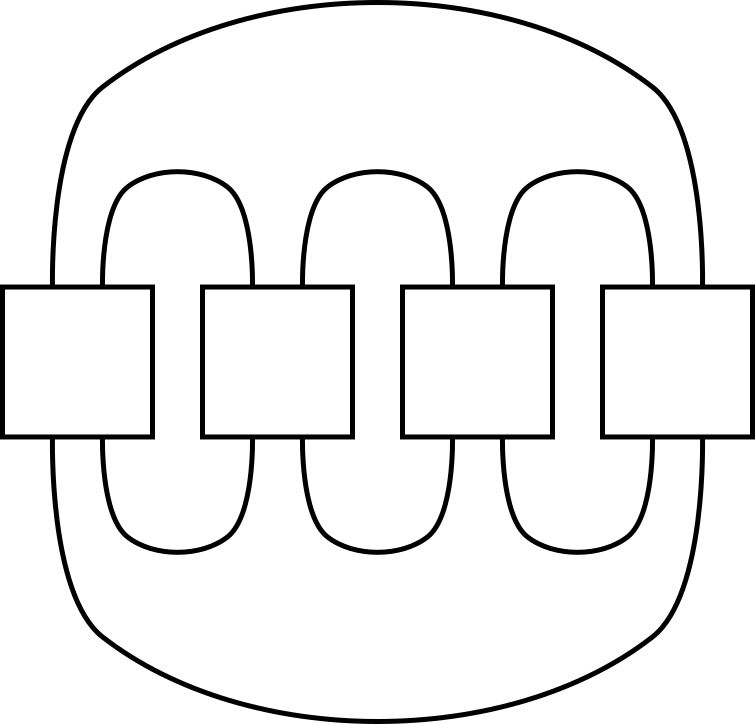}};

\node[label={$p$}] at (0.6,2.4){};
\node[label={$q$}] at (2.2,2.4){};
\node[label={$-p$}] at (3.8,2.4){};
\node[label={$-q$}] at (5.3,2.4){};

\end{tikzpicture}
    \caption{The butterfly link of $[K_p,K_q]$.}
    \label{butterfly}
\end{figure}

Boyle and Issa associate a directed strongly invertible knot $K$ with a $2$-componenent link $L_b(K)$, called the \emph{butterfly link} of $K$ (Definition 4.1 in \cite{boyle_issa}). The link $L_b(K)$ is obtained from $K$ by a band move along a band parallel to the chosen half-axis of $K$, in such a way that the linking number between the two components of $L_b(K)$ is zero.
Recall that a $n$-component link $L\subset S^3$ is said to be \emph{strongly slice} if it bounds $n$ disjoint disks properly embedded in $B^4$.
Boyle and Issa then prove the following result (see Proposition 7 in \cite{boyle_issa}).
\begin{prop}\label{prop}\footnote{This is a weak version of Proposition 7 in \cite{boyle_issa}. In particular Boyle and Issa define $L_b(K)$ to be a $2$-periodic link and prove that if $K$ is trivial in $\C$ then $L_b(K)$ is actually equivariantly slice.}
Let $K$ be a directed strongly invertible knot which is equivariantly slice. Then its butterfly link $L_b(K)$ is strongly slice.
\end{prop}
Hence, the proof of Theorem \ref{teorema} consists in showing that the butterfly link $L_b([K_p,K_q])$ is not strongly slice, which in turn relies on the following result of Aceto, Kim, Park and Ray.
\begin{thm}[\cite{aceto_kim_park_ray_2021}, Theorem 1.2]\label{pretzel}
Let $p$ and $q$ be odd integers such that $1<p<q$. Then the $2$-components pretzel link $P(p,q,-p,-q)$ is not strongly slice.
\end{thm}
Observe that $K_p^{-1}$, represented in Figure \ref{inverse}, is equivariantly isotopic to $K_{-p}$ through the $\pi$-rotation around the vertical axis in the figure.
By performing the sequence of equivariant connected sums, we see that the directed strongly invertible knot in Figure \ref{commutator} (where the chosen half-axis is the solid one) represents the commutator $[K_p,K_q]=K_p\widetilde{\#}K_q\widetilde{\#}K_p^{-1}\widetilde{\#}K_q^{-1}$.
We obtain now the butterfly link $L_b([K_p,K_q])$ in Figure \ref{butterfly} by a band move along a band parallel to the chosen half-axis of $[K_p,K_q]$ (the solid one in Figure \ref{commutator}).
Since $L_b([K_p,K_q])$ is the pretzel link $P(p,q,-p,-q)$, by Theorem \ref{pretzel} we know that it is not strongly slice.
Therefore by Proposition \ref{prop} the commutator $[K_p,K_q]$ is not equivariantly slice. This concludes the proof of Theorem \ref{teorema}.\hfill$\square$
\medskip\par

\section{Further remarks}
\subsection{An independent proof}\label{remark}
\begin{figure}[ht]
\centering
\begin{tikzpicture}

\node[anchor=south west,inner sep=0] at (0,0){\includegraphics[scale=0.3]{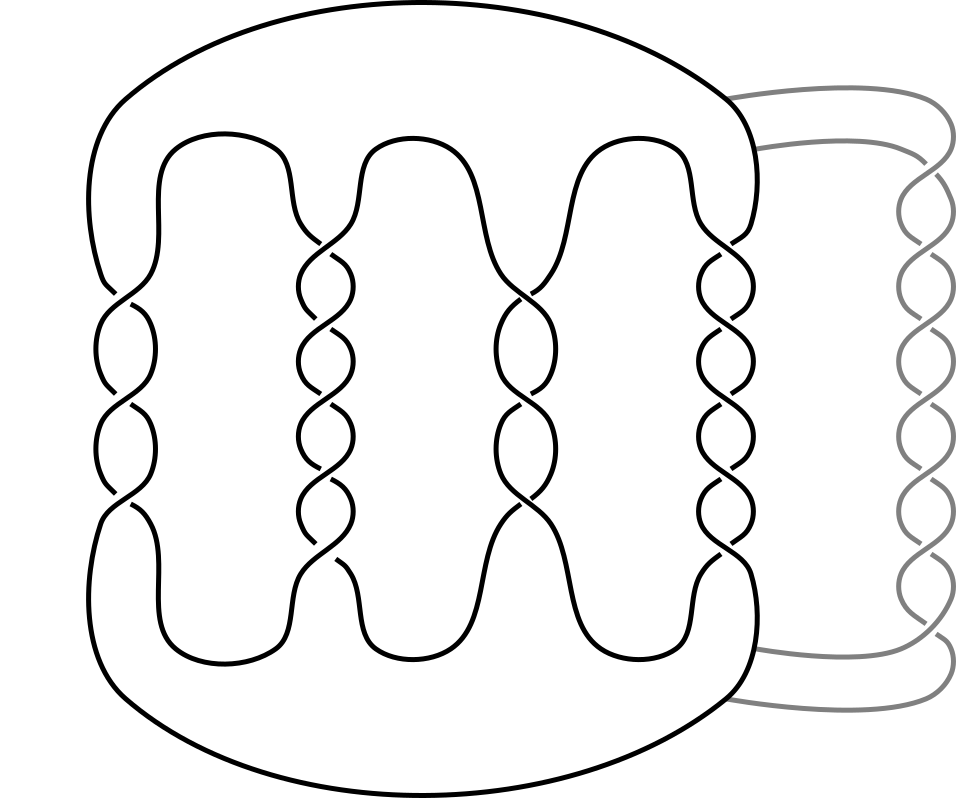}};
\end{tikzpicture}
    \caption{The pretzel link $P(3,5,-3,-5)$ and the band (in grey) that gives $P(3,5,-3,-5,7)$.}
    \label{pretzel35}
\end{figure}

Herald, Kirk and Livingston (Section 11 in \cite{Herald2008MetabelianRT}) prove that the pretzel knot $P(3,5,-3,-5,7)$ is not topologically slice by use of the twisted Alexander polynomials. This result leads to a proof that the pretzel link $P(3,5,-3,-5)$ is not topologically strongly slice, and in turn to an independent proof that $[K_3,K_5]$ is nontrivial in $\C$, by the same argument used in Section \ref{proof}. 
In fact, observe that $P(3,5,-3,-5,7)$ is obtained by a band move on $P(3,5,-3,-5)$ which connects the two components of the link, as in Figure \ref{pretzel35}.
This band move can be seen as a cobordism $C\subset S^3\times[0,1]$ of genus $0$ (i.e. a pair of pants) between $P(3,5,-3,-5)$ and $P(3,5,-3,-5,7)$. If $P(3,5,-3,-5)$ bounds a pair of (locally flat) disjoint disks $D_1\sqcup D_2$ in $B^4$ then $D=(D_1\sqcup D_2)\cup C\subset B^4\cup (S^3\times[0,1])\cong B^4$ would be a topological slice disk for $P(3,5,-3,-5,7)$, in contradiction with the result from \cite{Herald2008MetabelianRT}.

\subsection{The theta-curves cobordism group is not abelian}
Sakuma in \cite{sakuma} uses the relation between strongly invertible knots and $\theta$-curves in order to prove, in particular, that the equivariant connected sum of directed strongly invertible knot is not commutative.
The concordance of $\theta$-curves (in the piecewise linear category) was studied as independent topic by Tanyama in \cite{taniyama_1993}, who defined the \emph{cobordism group of $\theta$-curves} $\Theta$.
Miyazaki in \cite{miyazaki} provided a proof that the group $\Theta$ is not commutative, appealing on a result of Gilmer \cite{gilmer}. However, Friedl \cite{friedl} found gaps in the proof of the result in \cite{gilmer}.

We observe that the example in Subsection \ref{remark} of a nontrivial commutator in $\C$ can be adapted to provide a different proof that $\Theta$ is not abelian.
Notice that the union of a directed strongly invertible knot $K$ with its chosen half-axis produces a $\theta$-curve $\theta(K)$, with vertices ordered as initial and final point of the chosen half-axis (as in Section \ref{intro}).
It is easy to check that this defines a homomorphism
$$\theta:\C\longrightarrow\Theta.$$
Tanyama proved that a $\theta$-curve is slice if and only if any (and hence all) of its \emph{parallel link} is strongly slice (see Theorem 5 in \cite{taniyama_1993}).
Given a directed strongly invertible knot $K$, one of the parallel link of $\theta(K)$ is easily seen to be exactly the butterfly link $L_b(K)$.
Therefore, using the result from \cite{Herald2008MetabelianRT} described in Subsection \ref{remark}, we get that the commutator $[\theta(K_3),\theta(K_5)]$ is nontrivial in $\Theta$, and hence that the cobordism group of $\theta$-curves is not abelian.
Note that in this case \cite{aceto_kim_park_ray_2021} cannot be applied because its results only hold in the smooth category.
\vspace{4mm}

\paragraph*{\bfseries Acknowledgments} I would like to thank my friends and colleagues who read the drafts of this work for their helpful comments, and Paolo Lisca for the help and support during my master thesis project, which led to this paper.
I also thank Keegan Boyle for his helpful comments.

\bibliographystyle{amsalpha}

\bibliography{main}

\providecommand{\bysame}{\leavevmode\hbox to3em{\hrulefill}\thinspace}
\providecommand{\MR}{\relax\ifhmode\unskip\space\fi MR }
\providecommand{\MRhref}[2]{%
  \href{http://www.ams.org/mathscinet-getitem?mr=#1}{#2}
}
\providecommand{\href}[2]{#2}
\begin{thebibliography}{AKPR21}

\bibitem[AB21]{alfieri_boyle}
Antonio Alfieri and Keegan Boyle, \emph{Strongly invertible knots, invariant
  surfaces, and the {Atiyah-Singer} signature theorem}, 2021.

\bibitem[AKPR21]{aceto_kim_park_ray_2021}
Paolo Aceto, Min~Hoon Kim, Junghwan Park, and Arunima Ray, \emph{Pretzel links,
  mutation, and the slice-ribbon conjecture}, Mathematical Research Letters
  \textbf{28} (2021), no.~4, 945--966.

\bibitem[BI21]{boyle_issa}
Keegan Boyle and Ahmad Issa, \emph{Equivariant 4-genera of strongly invertible
  and periodic knots}, 2021.

\bibitem[DMS22]{dai_mallick_stoffregen}
Irving Dai, Abhishek Mallick, and Matthew Stoffregen, \emph{Equivariant knots
  and knot {Floer} homology}, 2022.

\bibitem[Fri04]{friedl}
Stefan Friedl, \emph{{Eta invariants as sliceness obstructions and their
  relation to Casson–Gordon invariants}}, Algebraic \& Geometric Topology
  \textbf{4} (2004), no.~2, 893 -- 934.

\bibitem[Gil83]{gilmer}
Patrick~M. Gilmer, \emph{{Slice knots in $S^3$}}, The Quarterly Journal of
  Mathematics \textbf{34} (1983), no.~3, 305--322.

\bibitem[HKL08]{Herald2008MetabelianRT}
Christopher Herald, Paul~A. Kirk, and Charles Livingston, \emph{Metabelian
  representations, twisted {Alexander} polynomials, knot slicing, and
  mutation}, Mathematische Zeitschrift \textbf{265} (2008), 925--949.

\bibitem[Miy94]{miyazaki}
Katura Miyazaki, \emph{{The Theta-Curve Cobordism Group Is Not Abelian}}, Tokyo
  Journal of Mathematics \textbf{17} (1994), no.~1, 165 -- 169.

\bibitem[Sak86]{sakuma}
Makoto Sakuma, \emph{On strongly invertible knots}, Algebraic and topological
  theories (Kinosaki, 1984) (1986), 176--196.

\bibitem[Tan93]{taniyama_1993}
Kouki Taniyama, \emph{{Cobordism of theta curves in $S^3$}}, Mathematical
  Proceedings of the Cambridge Philosophical Society \textbf{113} (1993),
  no.~1, 97–106.

\end{thebibliography}

\end{document}